\documentstyle[12pt]{article}

\begin{document}

\begin{center}
{\large \bf Gradient K\"ahler-Ricci 
Solitons and Periodic Orbits}
\end{center}
\begin{center}
\begin{tabular}{ccc}
Huai-Dong Cao \& Richard S. Hamilton\\
\end{tabular}
\end{center}
\footnotetext[2]{The authors' research was supported in part by NSF. }
\begin{abstract} 
We study Hamiltonian dynamics of gradient K\"ahler-Ricci solitons
that arise as limits of dilations of singularities of the Ricci flow on compact
K\"ahler manifolds. Our main result is that the underlying spaces of such 
gradient solitons must be Stein manifolds. Moreover, on all most all energy 
surfaces of the potential function $f$ of such a soliton,  the Hamiltonian vector
field $V_f$ of $f$, with respect to the K\"ahler form of the gradient soliton 
metric, admits a periodic orbit. 
The latter  should be of significance in the study of singularities of the Ricci flow on compact K\"ahler 
manifolds in light
of the ``little loop lemma'' principle in \cite{Has}. 
\end{abstract}

\vspace{1 cm}

\section{Introduction}

Let $X$ be a noncompact complex manifold. Consider the Ricci flow 
\begin{eqnarray}
{\partial \over {\partial t}}g_{i\bar j} =- R_{i\bar j}
\label{1}
\end{eqnarray}
which evolves a complete K\"ahler metric $g_{i\bar j}$ on $X$ by its 
Ricci tensor $R_{i\bar j}$.
A solution metric $g_{i\bar j}$ of Eq.(1) is called a gradient K\"ahler-Ricci 
soliton if it moves along 
Eq.(1) under a one-parameter family of biholomorphisms of $X$ generated 
by a holomorphic gradient vector field. Namely, 
$g_{i\bar j}(t)={\phi}^{*}(t)g_{i\bar j}(0)$, where $\phi (t)$
is a $1$-parameter family of automorphisms of $X$ generated by some 
holomorphic gradient vector field. More precisely, this means that 
the  Ricci tensor of the metric $g_{i\bar j}$ can be expressed as 
\begin{eqnarray}  
 R_{i\bar j}=f_{i\bar j}, \qquad \mbox {with} \qquad f_{,ij}=0  
\label{2}
\end{eqnarray}
for some real-valued function $f$ on $X$.
Note that the condition $f_{,i j}=0$ is equivalent to saying that 
the gradient vector field $\nabla f$ is holomorphic. 

Understanding gradient Ricci solitons is very important in the study of 
Ricci flow, since they often arise as limits of dilations of singularities of 
the Ricci flow. This was first observed by the second author \cite{Has}.
For example, consider the normalized Ricci flow
\begin{eqnarray}
{\partial \over {\partial t}}g_{i\bar j} =- R_{i\bar j}+g_{i\bar j}
\label{3}
\end{eqnarray}
on a $n$-dimensional compact K\"ahler manifold $M^n$ with positive first Chern 
class $C_1(M)$. It is known (\cite{Ca85}, \cite{Ca92}) that the solution to Eq.(3) exists for all time 
$0\leq t < \infty$. A key problem in the study of asymptotic behavior of solutions to Eq.(3) 
is to obtain a uniform estimate of the evolved curvature tensor. By singularities of the Ricci flow (3) we mean solutions 
to Eq.(3) with unbounded curvature as time $t$ tends to the infinity.
It is conjectured that singularities of Eq.(3) do not occur and solutions to Eq.(3)
either converges to a K\"ahler-Einstein metric or a homothetically shrinking 
gradient K\"ahler-Ricci soliton of Eq.(3), a solution moving along Eq.(3) under a  one-parameter family of biholomorphisms of $X$ generated 
by a holomorphic gradient vector field.
On the other hand, if singularities do 
occur  then, as indicated in  \cite{Has}, certain limit of parabolic 
dilations, dilating time as distance squared, of the solution metric 
$g_{i\bar j}(t)$ to Eq.(3) is a complete solution metric to Eq.(1) on a noncompact complex manifold $X$ 
defined for $-\infty <t<\infty$ with uniformly bounded curvature. When the 
holomorphic bisectional curvature of the initial metric
is positive this limit metric has nonnegative holomorphic 
bisectional curvature where the scalar curvature $R$ assumes its maximum in 
space-time.
It is shown in \cite{Ca97} that such a solution metric $g_{i\bar j}$ is necessarily a gradient 
K\"ahler-Ricci soliton. 
Thus singularities of the Ricci flow (3) are related to these special gradient K\"ahler-Ricci solitons.
So far, the only known examples of gradient 
K\"ahler-Ricci solitons are the rotationally symmetric ones.  They are
the "cigar" soliton on the complex plane found by the second author \cite{Ha86} and its higher dimensional 
analogues on ${\bf C}^n$ found by the first author \cite{Ca94}. Not much was known about 
general structures of gradient solitons. 
In this paper, we explore the structure of gradient solitons and the symplectic geometry of periodic orbits of certain
special Hamiltonian functions on gradient 
K\"ahler-Ricci solitons, which arise as limits of dilations of singularities of 
the Ricci flow (3). Our main result can be stated as follows. \\

\noindent {\bf Theorem 1.1} {\em Let $g_{i\bar j}$ be a complete gradient 
K\"ahler-Ricci soliton on a noncompact complex manifold $X$ with positive Ricci curvature such that the scalar curvature
$R$ assumes its maximum in space-time. Let $f$ be the 
potential function of the soliton as defined in (\ref{2}). Then $f$ is a plurisubharmonic 
exhaustion function. Hence $X$ must be a Stein manifold. Furthermore, near each level surface $S_c=\{f=c\}\subset X$ there
exists a periodic orbit for the Hamiltonian vector field $V_f$ of $f$ with 
respect to the symplectic form $\omega$ defined by the K\"ahler form of 
$g_{i\bar j}$.}\\

Our work is motivated by the second author's ``little loop lemma'' principle 
\cite{Has}, which roughly says that there is no small loop in a big flat 
region for the solution of Ricci flow whenever Hanarck estimates 
(\cite{Hah}, \cite{Ca92}) hold. This turns out to be an effective tool to rule
out singularities in some cases. Our work is also
motivated by the examples of the cigar soliton and its higher dimensional 
versions, where periodic orbits of the corresponding Hamiltonian vector fields 
do exist and have uniformly finite length. 

The existence of periodic orbits on a compact energy surface is an old question
of the qualitative theory of Hamiltonian systems. We refer the reader to the book
of Hofer and Zehnder \cite{HZ} for history and recent development on this 
rich subject. 

The proof of Theorem 1.1 is based on the structure equations of gradient 
Ricci solitons, the symplectic capacity theory of Hofer-Zehnder
(see \cite{HZ}), and a result of Eliashberg-Gromov \cite{EG}. One interesting 
feature is that we use both the K\"ahler form $\omega$ and the 
Ricci form $\rho=\sqrt{-1}\partial\bar\partial f$ of the soliton metric as 
symplectic forms on $X$. On one hand,
we need to use $\omega$ as the symplectic form due to geometric considerations.
On the other hand, the Ricci form $\rho$ is more naturally adapted to
symplectic embedding problems.  It turns out that $X$, with the Ricci form 
$\rho$ as the 
symplectic form,
can be symplectically embedded into $I\!\! R^{2n}$ with the standard
symplectic form.  Moreover, based on a result of Eliashberg-Gromov in \cite{EG},
we are able to show that there exists a symplectic embedding of a neighborhood
of the level 
surface $S_c$, with symplectic form $\omega$, into the symplectic manifold
$(X, \rho)$. It then follows from the capacity theory of Hofer-Zehnder \cite{HZ} 
that level surface $S_c$ has a periodic orbits for almost all $c>0$. 

\medskip

{\em Acknowledgement} The first author would like to thank
Yasha Eliashberg for valuable suggestions and Viktor Ginzburg for helpful 
discussions. 

\section {Gradient K\"ahler-Ricci solitons}

In this section we discuss various aspects of gradient K\"ahler-Ricci solitons.
First, we derive the first order and the second order equations 
satisfied by gradient Ricci solitons. To certain extent, these equations characterize
gradient K\"ahler-Ricci solitons (see \cite{Ca97}). Next, we describe the rotatinally symmetric
examples, namely, the cigar soliton on ${\bf C}$ and its analogues on ${\bf C}^n$.

\subsection {Ricci soliton equations}

Let $g_{i\bar j}$ be a gradient K\"ahler-Ricci soliton of Eq.(1)
so that it satisfies Eq.(2). 
From the commutation formulas for the covariant derivatives we  have
\[f_{,ik\bar j}-f_{,i\bar j k}=R_{i\bar j k\bar q}f_{,q}.\]
On the other hand,  taking covariant derivate of the first equation in (\ref{2}) 
and using the second equation,
we have 
\[f_{,ik\bar j}-f_{,i\bar j k}=-R_{i\bar j, k}.\]
Hence we get
\begin{eqnarray}  
R_{i\bar j, k}+R_{i\bar j k\bar q}f_{,q}=0.
\label{4}
\end{eqnarray}
Taking the trace, we get 
\begin{eqnarray}  
R_{, k}+R_{ k\bar l}f_{,l}=0.
\label{5}
\end{eqnarray}

Differentiating (\ref{4}) once more we obtain
\[R_{i\bar j, k\bar l}+R_{i\bar j k\bar q,\bar l}f_{,q}
+R_{i\bar j k\bar q,}f_{,q\bar l}=0.\]
By the  second Bianchi
identity and (\ref{2}) we have
\begin{eqnarray*}
R_{i\bar j, k\bar l}+R_{i\bar j k\bar l,\bar q}f_{,q
}+R_{i\bar j k\bar q}R_{q\bar l}=0.
\end{eqnarray*}
Taking the trace in the above equation, we obtain
\begin{eqnarray}
\Delta R_{i\bar j}
+R_{i\bar j k\bar l}R_{l\bar k}+R_{i\bar j,\bar k}f_{,k}=0.
\label{6}
\end{eqnarray}

Equations (\ref{4}) and (\ref{6}) are the first order and the second order
equations satisfied by gradient K\"ahler-Ricci solitons. They are also related
to the (1,1) tensor in the Hanarck estimates (see \cite{Ca92}).

\subsection {Rotationally symmetric examples}

In \cite{Ha86}, the second author wrote down the first example of a gradient Ricci 
soliton, called the cigar soliton, on the complex plane ${\bf C}$. It has 
the form
\[ds^2={|dz|^2\over {1+|z|^2}},\qquad z\in {\bf C}. \]
Its potential function $f$ is given by $f=\log (1+|z|^2)$ and its scalar 
curvature  is given by 
\[R=\frac {1}{1+|z|^2}.\]
 It is easy to verify that $\nabla f=z$. Here $\nabla f$ denotes the gradient
 of $f$ with respect to the soliton metric $ds^2$. So the cigar soliton flows towards the origin 
by conformal dilations and has positive 
Gaussian curvature.
It is asymptotic to a flat cylinder at infinity and has maximal scalar curvature
$1$ at the origin, hence the name cigar soliton. 

In \cite{Ca94}, the first author studied the existence of rotationally symmetric gradient
K\"ahler-Ricci solitons and proved the following\\

\noindent{\bf Proposition 2.1} {\em For  $n\ge 1$, there exists on 
${\bf C}^{n}$
a complete rotationally symmetric gradient K\"ahler-Ricci 
soliton of positive sectional curvature.}\\

Our rotationally symmetric solitons are constructed by considering 
$U(n)$-invariant K\"ahler metrics of the form 
\[g_{i\bar j}={\partial}_i{\partial}_{\bar j}\Phi (z,\bar z), \]
with the K\"ahler 
potential $\Phi (z,\bar z)$ expressed as
\[ \Phi (z,\bar z) =u(t), \qquad \bigl(t=\log |z|^2\bigr).\] 
Here $u(t)$ is a smooth function  defined on $(-\infty, \infty)$ 
satisfying the differential inequalities
\begin{eqnarray}
 u'(t)>0,\quad u''(t)>0, \qquad t\in (-\infty, \infty)
\label{}
\end{eqnarray} 
and as 
 $t\rightarrow -\infty$, has an expansion   
\begin{eqnarray}
u(t)=a_0+a_1e^t+a_2e^{2t}+\cdots. \qquad (a_1>0) 
\label{}
\end{eqnarray} 

Note that the metric $g_{i\bar j}$ is given by
\[g_{i\bar j}={\partial}_i{\partial}_{\bar j}u(t)=e^{-t}u'(t){\delta}_{i\bar j}
+e^{-2t}{\bar z}_iz_j\bigl(u''(t)-u'(t)\bigr).\]
It is natural to consider
\[f(t)=-\log \det (g_{i\bar j})=nt-(n-1)\log u'(t)-\log u''(t) \]
so that the Ricci tensor of the metric $g_{i\bar j} $ is given by
$R_{i\bar j}={\partial}_i{\partial}_{\bar j}f(t)$.

According to (2), the metric $g_{i\bar j}$ will be a gradient K\"ahler-Ricci soliton with 
potential function $f$ if the gradient vector field  
$V=\nabla f$ is holomorphic. By a direct computation, 
\[V^i=g^{i\bar j}f_{,\bar j}=z_i{f'(t)\over u''(t)}.\]
Hence $V$ is holomorphic if and only if 
$f'(t)=\alpha u''(t)$
for some constant $\alpha\in {I\!\! R}$. 
By setting $\phi (t)=u'(t)$ and by appropriate scaling (in metric) and 
dilation (in $z$), the equation satisfied by a gradient soliton 
becomes 
\begin{eqnarray}
{\phi}^{n-1}{\phi}'e^{\phi}=e^{nt}. 
\end{eqnarray}
An implicit solution of $\phi$ is given by
\begin{eqnarray}
\sum_{k=0}^{n-1} (-1)^{n-k-1}{n!\over k!}{\phi}^ke^{\phi}=e^{nt}+(-1)^{n-1}n!
\end{eqnarray}

It is easy to check from (10) that the solution $\phi=u'$  satisfies the required 
asymptotic condition (8) and the differential inequalities (7). 
Thus it gives rise to a $U(n)$-invariant gradient K\"ahler-Ricci soliton 
metric $g$ on ${\bf C}^n$. 
Furthermore, the sectional curvature of the soliton
metric turns out to be positive (see \cite{Ca94}). Note that when $n=1$, we have $\phi =\log (1+e^t)$ and the soliton metric is
nothing but the cigar soliton.

From (10) or (9),  it is  also easy to see that
\[\lim_{t\to \infty}t^{-1}\phi(t)=n, \qquad \lim_{t\to \infty}\phi'(t)=n\]
which in turn imply that the soliton metric $g$ is complete. Let $s$
denote the distance function from the origin with respect to $g$. Then
it can be shown that the soliton metric $g$ 
on  ${\bf C}^n$ ($n\geq 2$) satisfies the 
following two properties: (1) The volume of the geodesic ball $B_g(0,s)$ with respect to the metric $g$
grows like ${s}^n$; (2) The scalar curvature $R$ of the metric 
$g$ decays like $1/s$.

Finally, we point out that on this Ricci soliton
the sphere $S^{2n-1}$ at radius $s$ looks like an $S^1$ bundle over $CP^{n-1}$ 
where the $CP^{n-1}$ has diameter on the order of $\sqrt {s}$ while the $S^1$ 
fiber has diameter on the order of 1. This $S^1$ fiber is nothing but 
a periodic orbit for the Hamiltonian vector field $V_f=J\nabla f$, 
$J$ being 
the standard complex structure on ${\bf C}^{n}$, of 
$f$. For the cigar soliton these periodic orbits
are simply given by parallel circles, whose lengths are bounded by 
$2\pi r/\sqrt{1+r^2}$.
Here $r$ denotes the Euclidean radii of the circles.

\section {Symplectic capacity theory of Hofer-Zehnder}

In this section,  we recall the symplectic capacity function introduced by Hofer and Zehnder
(see \cite{HZ}) and the related existence theorem of periodic orbits on energy surfaces.
The results will be used in the next section to prove Theorem 1.1.

\subsection {Definition of the capacity $c_0$}

Consider the class of all symplectic manifolds $(M,\omega)$ possibly with boundary and of 
fixed dimension $2n$. A symplectic capacity is a map $(M,\omega) \rightarrow
c(M,\omega)$ which associates with every symplectic manifold $(M,\omega)$ a
nonnegative number or $\infty$, satisfying the following three properties.
\begin{itemize}
\item{(1)} Monotonicity: \ $c(M,\omega) \leq c(N,\tau)$ if there exists a symplectic embedding
$\phi :(M,\omega) \rightarrow (N,\tau)$.

\item{(2)} Conformality: \ $c(M,a\omega)=|a|(M,\omega)$ for all $a\in I\!\! R$,
$a\neq 0$.

\item{(3)} Nontriviality: \ $c(B_1,\omega)=\pi=(Z_1,\omega_0)$ for the open unit
$B_1$ and the open symplectic cylinder 
$Z_1=\{(x,y)\in I\!\! R^{2n}: x_1^2+y^2_1<1\}$ in the standard symplectic vector 
space $(I\!\! R^{2n}, \omega_0)$.
\end{itemize}

It follows from the conformality property and nontriviality property of a symplectic
capacity that for the open ball of radius $r>0$ in $I\!\! R^{2n}$, 
$c(B(r))=\pi r^2$. In particular, the capacity $c(K)$ of every compact subset $K$ in $I\!\! R^{2n}$ is finite.
 
Hofer and Zehnder \cite{HZ} introduced a distinguished capacity function $c_0$. The capacity 
$c_0(M, \omega)$ measures the minimal $C^{0}$-oscillation of special Hamiltonian 
functions 
$H: M\rightarrow I\!\! R$, needed in order to conclude the existence of a distinguished 
periodic orbit having small period. 

Now let's recall the definition of the capacity $c_0$.
On a symplectic manifold $(M,\omega)$, to a smooth function $H: M\rightarrow I\!\! R$
there belongs  a unique Hamiltonian vector field $V_H$ on $M$ defined by 
\[ (i_{V_{H}}\omega) (x)=\omega (V_H(x), w)=-dH(x)(w), \ \mbox {for} \ w\in T_x(M)\]
and $x\in M$. With respect to a given Riemannian metric $<,>$ on $M$ and 
a compatible almost complex structure $J$, the Hamiltonian vector field is
represented by 
\[ V_{H}(x)=J\nabla H(x)\in T_{x}M\]
where $\nabla H$ is the gradient of $H$ with respect to the metric $<,>$.
A $T$-periodic solution $x(t)$ of the Hamiltonian equation
\begin{eqnarray}
x'=V_H(x) \quad \mbox{on} \ M
\end{eqnarray}
is a solution satisfying the boundary conditions $x(T)=x(0)$ for some $T>0$.
Denote by $\mathcal{H}(M,\omega)$ the set of smooth functions $H$ on $M$
satisfying the following properties:

\begin{itemize}
\item{(a)} There is a compact set $K\subset (M\setminus\partial M)$ (depending on $H$) such 
that $H(x)\equiv m(H)$, a positive constant, on $M\setminus K$.

\item{(b)} There is an open set $U\subset M$ (depending on $H$) on which $H\equiv
0$.

\item{(c)} $0\leq H(x)\leq m(H)$ for all $x\in M$.

\item{(d)} All the periodic solutions of Eq. (11) on $M$ are either constant or 
have a period $T>1$. 
\end{itemize}
Hofer and Zehnder \cite{HZ} defined  
\[c_0(M,\omega)=\mbox{sup}\{m(H): H\in \mathcal{H}(M,\omega)\}\]
and proved \\
 
\noindent {\bf Theorem 3.1.} {\em The function $c_0$ is a symplectic capacity.}
\\

\subsection {A existence theorem of periodic orbits}

The flow $\phi^t$ of a Hamiltonian vector field
\[ x'=V_H(x), \ x\in M\]
leaves the level sets of the function $H$ on $M$ invariant. Fixing a value $c$ 
of this 
energy function which we can assume to be $c=1$, suppose that the energy surface 
\[ S=\{x\in M: H(x)=1\}\]
is compact and regular, i.e., $dH(x)\neq 0$, for $x\in S$. Then $S\subset M$ is a 
smooth and compact submanifold of codimension $1$ whose tangent space at $x\in S$ is given by 
\[T_xS=\{\xi\in T_xM: dH(x)\xi=0\}.\]
By the definition of a Hamiltonian vector field we have $V_H(x)\in T_xS$ so that $V_H$
is a nonvanishing vector field on $S$ whose flow exists for all time 
since $S$ is compact. On the other hand, since $S$ is compact and regular there
is an open and bounded neighborhood $U$ of $S$ which is filled with compact and regular 
energy surfaces having energy values near $c=1$:
\[U=\bigcup_{\lambda\in I} S_{\lambda},\]
where $I$ is an open interval around $\lambda =1$ and 
$S_{\lambda}=\{x\in U: H(x)=\lambda\}$ is diffeomorphic to the given surface $S=S_1$.

Now we can state a result of Hofer and Zehnder (\cite{HZ}, p. 106) about the existence of periodic 
orbits.\\ 


\noindent{\bf Theorem 3.2} {\em Let $S$ be a compact regular energy surface for the 
Hamiltonian vector field $V_H$ on $(M,\omega)$. Assume there is an open neighborhood
$U$ of $S$ having finite capacity: $c_0(U,\omega)<\infty$. 
Then there is a dense subset $\Sigma\subset I$ such that for $\lambda\in \Sigma$ the
energy surface $S_{\lambda}$ has a periodic solution of $V_{H}$.}

\section{The Proof of Theorem 1.1}

Let $g_{i\bar j}$ be a complete gradient K\"ahler-Ricci soliton on $X$ with 
positive Ricci curvature such that the scalar curvature
$R$ assumes its maximum $R_{max}=1$ at the point $O$, called the origin of $X$.
Let $f$ be the potential function of the gradient soliton $g_{i\bar j}$ 
defined in Eq.(\ref{2}).\\

\noindent{\bf Lemma 4.1} {\em Let $g_{i\bar j}$ be a complete gradient 
K\"ahler-Ricci soliton on $X$. Then we have 
\[
 R+|\nabla f|^2=C\]
on $X$, where $C$ is a positive constant.}\\

 {\bf Proof}: Differentiating $R+|\nabla f|^2$ and using (2) and (5), we obtain
 \[ (R+|\nabla f|^2)_i= R_{,i}+f_{j}f_{i\bar j}= R_{,i}+f_{j}R_{i\bar j}=0.\]
 Hence $R+|\nabla f|^2$ is a constant. 
 \hspace*{\fill}q.e.d. \\

\noindent {\bf Proposition 4.2} {\em The potential function $f$ is a strictly convex 
exhaustion function on $X$ and assumes its minimum at the origin $O$.
As a consequence, $X$ is necessarily a Stein manifold.}\\

{\bf Proof}. The argument is basically the same as in \cite{Has}. Since the Ricci curvature $R_{i\bar j}$ is positive, it follows
from the first equation in (2) that $f$ is plurisubharmonic. Furthermore, 
the second equation in (2) implies that $f$ is in fact strictly convex
on $X$. Thus it remains to show that $f$ is proper, i.e., the set $\{f\leq c\}$
is compact for every positive number $c>0$, and assumes its minimum at $O$. 
From Lemma 4.1 we know that
\begin{eqnarray}
R+|\nabla f|^2=C 
\end{eqnarray}
for some positive constant $C$. Hence $|\nabla f|$ is smallest at 
the origin $O$, since $R$ has maximum at $O$. We claim that $C=1$. If not, then $C>1$ and 
$|\nabla f|^2>0$ everywhere. Consider a gradient
path of $f$ through the origin and parametrized by $z^i=z^i(u)$ with $x^i$ at
the origin at $u=0$ and 
\[\frac {dz^i}{du}=g^{i\bar j}f_{\bar j}.\]
Then 
\[{d\over du}|\nabla f|^2=2g^{i\bar l}g^{k\bar j}R_{k\bar l}f_if_{\bar j}
> 0\]
since $R_{i\bar j}>0$ and $|\nabla f|^2>0$. Then $|\nabla f|^2$ isn't smallest at the 
origin, a contradiction. Therefore $C=1$  and hence $\nabla f=0$ at the origin.
Note also, from (12), that $|\nabla f|^2$ is bounded from 
above by $1$. 
These facts, together with the fact that $f$ is strictly convex, implies that $f$ is least at 
the origin and that $f$ is proper. Now a theorem of Grauert 
says  that $X$ must be a Stein manifold. 
\hspace*{\fill}q.e.d. \\

Note that a plurisubharmonic function $\phi$ on a Stein manifold $M$ gives rise
to a positive definite closed real (1,1) form 
\[\omega_{\phi}=\sqrt{-1}\partial\bar\partial\phi\],
and hence a symplectic form,  on $M$.\\

\noindent {\bf Lemma 4.3} {\em There exists a symplectic embedding 
$H :(X,\rho) \rightarrow ({\bf C}^{n},\omega_0)$}, where $\omega_0=
\frac{\sqrt{-1}}{2}\sum_i dz^i\wedge d{\bar z}^i$ is the standard symplectic 
form on ${\bf C}^n$, and $\rho=\sqrt{-1}\partial\bar\partial f$ is the Ricci form
of the soliton metric $g_{i \bar j}$.\\

{\bf Proof}. We may write $\rho=d\eta$ and $\omega_0=d\eta_0$, where $\eta$ 
and $\eta_0$ are real $1$-forms on $X$ and ${\bf C}^{n}$ respectively. From 
Darboux's theorem, there exists a symplectic diffeomorphism 
$h: (U,\rho) \rightarrow (U_0,\omega_0)$ from a neighborhood $U$ of $O\in X$
onto a neighborhood  $U_0$ of the origin in ${\bf C^{n}}$. We denote 
\[\tilde\eta_0=(h^{-1})^{*}(\eta)\]
the push-forward of $\eta$. Then the form $\eta_0-\tilde\eta_0$ is closed on $U_0$.
Hence $\eta_0-\tilde\eta_0=d\psi$ for some function $\psi$ in a smaller neighborhood $U_0'\subset U_0$.
Extend this function $\psi$, still denoted as $\psi$, to entire ${\bf C}^{n}$ 
so that $\psi\equiv 0$ outside $U_0$. Define the $1$-form 
\[ \eta'_0=\eta_0-d\psi.\] 
Then $\eta'_0=\eta_0$ outside $U_0$, $\eta'_0=\tilde\eta_0$ on $U_0'$ and 
\[d\eta'_0=d\eta_0=\omega_0.\]

Let $W'$ be the vector field, which is integrable, on ${\bf C}^n$ dual to 
$\eta'_0$ with respect to $\omega_0$. Then , by an easy argument, we can obtain a unique extension 
$H :(X,\rho) \rightarrow ({\bf C}^{n},\omega_0)$ of 
$h: (F^{-1}(U_0'),\rho) \rightarrow (U'_0,\omega_0)$ by following the trajectries of 
$W=\nabla_{g_{f}}f$ and $W'$ respectively. Here $\nabla_{g_{f}}f$ denotes the gradient
vector field of $f$ with respect the K\"ahler metric defined by the Ricci form
$\rho$.
\hspace*{\fill}q.e.d. \\

Next we need the following result of Eliashberg-Gromov in \cite{EG}.\\

\noindent{\bf Proposition 4.4} {\em Let $\phi$ and $\psi$ be two exhaustion
plurisubharmonic functions on a Stein manifold $M$. Then there exists a 
symplectic diffeomorphism $G :(M,\omega_{\phi}) \rightarrow (M,\omega_{\psi})$.}\\

\noindent {\bf Lemma 4.5} {\em For every compact subset $K\subset X$, 
there exists a symplectic embedding 
$G :(K,\omega) \rightarrow (X,\rho)$}.\\

{\bf Proof}. For simplicity, we assume that $K=\{f\leq c\}$ for some constant $c>0$.
Consider the compact set $K'=\{f\leq c'\}$ for some $c'>c$ so that $K\subset K'$.
We can write the K\"ahler form $\omega=\sqrt{-1}\partial\bar\partial\phi$ for some 
real-valued function $\phi$ on $X$, since $\omega$ is a closed (1,1) form  and $X$
is simply connected. We may assume $\phi\geq 0$
on $K$. Choose a convex diffeomorphism $h: {I\!\! R} \rightarrow {I\!\! R}$ such that
$\tilde f=h\circ f$ is again a plurisubharmonic exhaustion function on $X$, 
$\tilde f=f$ on $K$, and $\tilde f_{|{\partial K'}}$ very large. Define
\[\psi={\mbox max}(\tilde f-c, \phi)=\left\{
\begin{array}{ll}
   \phi & \mbox{on}\  K\\
 \tilde f-c &  \mbox{outside} \ K' 
\end{array}\right. \]
Then $\psi$ is another plurisubharmonic exhaustion function on $X$.
By Proposition 4.4, there exists a 
symplectic diffeomorphism $G :(X,\omega_{\psi}) \rightarrow (X,\rho=\omega_{f})$.
But on $K$ we have $\omega_{\psi}=\sqrt{-1}\partial\bar\partial\phi=\omega$.
Hence the restriction of $G$ to $K$ gives a symplectic embedding of 
$(K,\omega)$ into $(X,\rho)$.
\hspace*{\fill}q.e.d. \\

Combining Lemma 4.3 and Lemma 4.5, we obtain \\

\noindent{\bf Proposition 4.6} {\em Under the assumptions of Theorem 1.1, there exists
a symplectic embedding 
$F :(K,\omega) \rightarrow ({\bf C}^{n},\omega_0)$} (depending on K) for every compact 
subset $K\subset X$.\\



{\bf Proof of Theorem 1.1.} It follows from Proposition 4.2 that $X$ is a Stein 
manifold. From Proposition 4.6, we know that for every energy value $c>0$ we can pick a 
number $c'>c$ so that there exists a symplectic embedding 
$F :(K',\omega) \rightarrow ({\bf C}^{n},\omega_0)$ 
for the compact subset $K'=\{f\leq c'\}\subset X$. Since every compact subset 
of ${\bf C}^n$ has finite capacity, the monotonicity property of the  
capacity function $c_0$ implies that $K'$ has 
finite capacity $c_0(K',\omega)$ hence Theorem 3.2 applies. This completes the proof
of Theorem 1.1. 
\hspace*{\fill}q.e.d. \\

\noindent Texas A\&M University, College Station \\
{\em E-mail address}:  cao@math.tamu.edu\\

\noindent University of California, La Jolla and Columbia University, New York

\end{document}